# Randomized Algorithms for Solving Singular Value Decomposition Problems with MATLAB Toolbox

©2021
Xiaowen Li

Submitted to the graduate degree program in Department of Mathematics and the Graduate Faculty of the University of Kansas in partial fulfillment of the requirements for the degree of Master of Art.

# Abstract

This thesis gives an overview of the state-of-the-art randomized linear algebra algorithms [3, 7, 9, 11, 13, 14, 15, 18] for singular value decomposition (SVD), including presentation of existing pseudo-codes and theoretical error analysis. Our main focus is on presenting numerical experiments illustrating image restoration using various randomized singular value decomposition (RSVD) methods; theoretical error bounds, computed errors, and canonical angles analysis for these RSVD algorithms.

This thesis also comes with a newly developed MATLAB toolbox that contains implementations and test examples for some of the state-of-the-art randomized numerical linear algebra algorithms introduced in [3, 7, 9, 11, 13, 14, 15, 18].

# Acknowledgements

I would like to thank my master's research advisor – Dr. Agnieszka Międlar for her invaluable supervision, support and tutelage during the course of my Master's degree. Additionally, I would like to express gratitude to my undergrad advisor Dr. Ron Morgan from Baylor University to lead me to the field of numerical analysis. My appreciation also goes out to my family and friends for their encouragement and support all through my studies especially during coronavirus time periods.

# Contents

# List of Figures

# List of Algorithms

# Chapter 1

# Introduction

## 1.1 Standard Numerical Linear Algebra

Numerical Linear Algebra (NLA) plays an important role in applied mathematics, statistics and computer science. Defined as

> *The study of algorithms for the problems of continuous mathematics.*
> - L. N. Trefethen, Oxford University (1992)[16],

numerical linear algebra, in general, involves studies of the following problems:

- **Eigenvalue problems:** Given a matrix $A \in \mathbb{R}^{n \times n}$ find a scalar $\lambda \in \mathbb{R}$ (eigenvalue) and a nonzero vector $x \in \mathbb{R}^n$ (eigenvector) such that $Ax = \lambda x$.

- **Linear Systems:** Given a matrix $A \in \mathbb{R}^{n \times n}$ and vector $b \in \mathbb{R}^n$ find a vector $x \in \mathbb{R}^n$ such that $Ax = b$.

- **Least-Squares Problem:** Given a matrix $A \in \mathbb{R}^{m \times n}$ and vector $b \in \mathbb{R}^m$ find a vector $x \in \mathbb{R}^n$ such that $\min_x \|b - Ax\|_2$.

Term *Numerical Linear Algebra* was first introduced in a 1947 paper by von Neumann and Go [17]. Standard NLA algorithms such as Krylov subspace iteration methods or popular matrix factorization methods are developed and fully reliable for small to medium scale input matrix problems. Modern applications, e.g., machine learning and modern massive data set (MMDS) analysis, create challenges with much larger data sets and sometimes missing information. Randomization seems to be a promising technique to address these challenges.

## 1.2 Randomized Numerical Linear Algebra

Randomized Numerical Linear Algebra (RNLA) is an interdisciplinary study field that lets randomization into the classical already existing algorithms and adapt them to handle more efficiently large-scale linear algebra problems. From a foundational perspective, RNLA has its roots in the field of theoretical computer science, with deep connections to mathematics (convex analysis, probability theory, or metric embedding theory) and applied mathematics (scientific computing, signal processing, and numerical linear algebra) [3]. The efficiency relates to two major concerns: run-time and storage. In the real-world, the clock speed is hardly improving for hardware, but with RNLA applied, we are able to keep the run-time down from the algorithm side; the cost of slow storage (hard drives, flash memory, etc.) is also going down rapidly because of the RNLA algorithms. Moreover, RNLA provides a sound algorithmic and statistical foundation for cloud computing and big data analysis.

RNLA addresses a variety of topics such as randomized matrix multiplication [2], randomized least-squares solvers [5], and low-rank approximation [10], to mention just a few. This thesis will mainly focus on introducing randomized algorithms for singular value decomposition with low-rank matrices, with a MATLAB toolbox that contains implementations and test examples for basic randomized algorithms.

# Chapter 2

# Background/ Preliminary

## 2.1 Randomization

The reason that randomization is involved is because we want to take the input matrix and construct a `sketch` of this matrix with random sampling. Such sketch is a smaller or sparser matrix that can still represent the essential information contained in the original input matrix. Therefore all later computations can cost less but still give a reasonable accurate solution.

Such randomization has been considered and commonly used after the following observation:

**Johnson-Lindenstrauss Lemma**  From a 1984 paper [8], the Johnson-Lindenstrauss Lemma states that any $n$ points in high dimensional Euclidean space can be mapped onto $k$ dimensions where $k \geq \mathcal{O}(log n/\epsilon^2)$ without distorting the Euclidean distance between any two points more than a factor of $1 \pm \epsilon$. [1]

**Lemma 1** *(Johnson-Lindenstrauss [1]) Let $\epsilon$ be a real number such that $\epsilon \in (0,1)$ , $n$ be a positive integer and $k$ be an integer such that*

$$k \geq 4\Big(\frac{\epsilon^2}{2} - \frac{\epsilon^3}{3}\Big)^{-1} log(n).$$

*Then for any set $V$ of $n$ points in $\mathbb{R}^d$, there is a map $f : \mathbb{R}^d \to \mathbb{R}^k$ such that*

$$(1-\epsilon)\|u-v\|^2 \leq \|f(u)-f(v)\| \leq (1+\epsilon)\|u+v\|^2, \qquad \forall u, v \in V$$

*Further this map can be found in randomized polynomial time. Repeating this projection $\mathcal{O}(n)$ times can boost the success probability to any desired constant, giving us a randomized polynomial time algorithm.*

Lemma 1 shows that random embedding preserves Euclidean geometry. Which suggests that we should be able to solve many computational problems of a geometric nature more efficiently by translating them into a lower- dimensional space using sketching.

### 2.1.1 Types of randomization

There are several examples of considering randomization.

**Element-wise sampling (Sparsification)**   We can review an $m \times n$ matrix $A$ as a 2D array that contains $mn$ elements. We can denote each of such elements by $a_{ij}$ where $i = 1, ..., m; j = 1, ...n$. Then to in order to produce a `sketch` that is less costly for computations comparing to the original matrix size, instead of keeping all its elements, we can randomly select a portion of the elements, and let the rest of the matrix to be 0. Again, the sketch should be smaller but should still contain enough information, hence the way of sampling the entries is important. One simple and common choice is to pick entries with probability proportional to their squared-magnitudes.

**Row/column sampling**   Every element in a matrix can be specified by it's row and column index. Therefore, we can also randomly pick a portion of the matrix by selecting rows and columns. This method leads to much stronger worst-case bounds and is more commonly used in randomized linear algebra [4] due to its lower complexity.

## 2.2   Two-Stage Approach

Given an $m \times n$ input matrix A, a target rank $1 \le k \le rank(A)$, and an oversampling parameter $p \ge 0$. Our goal is to use randomization to compute a low-rank, e.g., rank $k + p$ approximation of matrix *A*, which can be performed in two stages:

**Stage 1:** Here, we obtain a low-dimensional subspace which contains the most important information in matrix *A*. We start with generating an random $n \times (k+p)$ Gaussian test matrix $\Omega$. Then we form a matrix product $Y = A\Omega$. Since $Y$ is likely to be ill-conditioned, we orthonormalize its columns to form an orthonormal basis $Q$ of the low-dimensional space of interest. This procedure is illustrate in Algorithm 1.

---
**Algorithm 1** Solving the Fixed Rank Problem[7]
---
    *Inputs:* An $m \times n$ matrix *A*, a target rank *k*, an oversampling parameter *p*.
    *Outputs:* An $m \times (k+p)$ orthonormal matrix *Q* which approximates the range of *A*
1: Generate an $n \times (k+p)$ Gaussian test matrix $\Omega$.
2: Form $Y = A\Omega$.
3: Construct a matrix *Q* whose columns form an orthonormal basis for the range of *Y*.
---

**Stage 2**: In this stage we take advantage of the lower-dimensional space determined in Stage 1. We first obtain a reduced matrix by restricting an original matrix *A* to the lower-dimensional space and then compute its desired factorization, e.g., the singular value decomposition (SVD) (see Chapter 3).

# Chapter 3

# Randomized Singular Value Decomposition

In this chapter we will discuss several different variants of determining the *randomized singular value decomposition.*

## 3.1 General Randomized SVD Algorithm

In Section 2.2, we have introduced a general two-stage procedure to obtain a low-rank factorization of a given matrix $A$ using randomization. Let us illustrate this approach on the example of singular value decomposition. In particular, let us concentrate on describing the Stage 2. In the case of general randomized SVD method, we start with forming matrix $B = Q^*A$. Since matrix $Q$ computed in Stage 1 satisfy $\|A - QQ^*A\|_2 \leq \varepsilon$, where $\varepsilon$ denotes a small number that is less than 1, $A \approx QB$. Then we apply the standard algorithm to obtain the SVD of a small matrix $B$, i.e., $B = \widetilde{U}\Sigma V^*$. Finally, by setting $U = Q\widetilde{U}$ we complete Stage 2 with $A \approx U\Sigma V^*$.

Given the $j$th largest singular value $\sigma_j$ of matrix $A$, we know that for each $j \geq 0$ the minimizer $X$ satisfies

$$\min_{rank(X) \leq j} \|A - X\|_2 = \sigma_{j+1}. \tag{3.1}$$

**Theorem 2** *[7, Theorem 1.1] Suppose that A is a real m × n matrix. Select a target rank $k \geq 2$ and an oversampling parameter $p \geq 2$, where $k + p \leq \min\{m, n\}$. Execute Algorithm 1 with a standard Gaussian random test matrix $\Omega$ to obtain an m × (k + p) matrix Q with orthonormal columns. Then*

$$\mathbb{E}\|A - QQ^*A\| \leq \left[1 + \frac{4\sqrt{k+p}}{p-1} \cdot \sqrt{\min\{m,n\}}\right]\sigma_{k+1},$$

*where* $\mathbb{E}$ *denotes the expectation with respect to the random test matrix and* $\sigma_{k+1}$ *is the (k+1)th largest singular value of the matrix A.*

From Theorem 2 we know that the Stage 1 of our randomized SVD algorithm (1) generates matrix $Q$ such that the resulting approximation of $A^*A$ is within a polynomial factor $\left(\left[1 + \frac{4\sqrt{(k+p)\min\{m,n\}}}{p-1}\right]\right)$ of the theoretical minimum. So if we fix the size of matrix $A$ as well as the target rank $k$, the error is only determined by the singular values of matrix $A$.

---

**Algorithm 2** General Randomized SVD Algorithm [13]

*Inputs:* An $m \times n$ matrix $A$, a target rank $k$, an oversampling parameter $p$.

*Outputs:* Orthonormal matrices $U$, $V$, and a diagonal matrix $\Sigma$ from a rank-$(k+p)$ SVD approximation of $A$.

**Stage 1:**
1: Generate an $n \times (k+p)$ Gaussian test matrix $\Omega$.
2: Form $Y = A\Omega$.
3: Construct a matrix $Q$ whose columns form an orthonormal basis for the range of $Y$.

**Stage 2:**
4: Form $B = Q^*A$.
5: Compute an SVD of the small matrix: $B = \widetilde{U}\Sigma V^*$.
6: Set $U = Q\widetilde{U}$.

---

If a factorization with exactly rank $k$ is desired, we can truncate the last $p$ components of the SVD factorization of matrix $B$ after line 5 in Algorithm 2. Larger oversampling parameter $p$ provides a more accurate output, but at the same time also less efficient. So, we can adjust the value of $p$ to trade accuracy with efficiency.

### 3.1.1 Complexity

To analyze the cost of the general randomized SVD method presented in Algorithm 2, there are several terms that we need to define. First, depending on different sampling techniques used in implementation, the exact cost for the sampling part of the algorithm varies. In general, the cost of this step is $\mathcal{O}(n(k+p))$ [6].

Second, the cost of matrix multiplication can also vary depending on the structure of the matrices. Therefore, we denote the cost for a matrix vector multiplication with matrix $A$ as $\tau_A$ [6]. Then we are able to analyze the cost of the general randomized SVD algorithm as follows:

- $\tau_A(k+p)$ cost for a $A\Omega$ matrix multiplication.

- $\mathcal{O}(m(k+p)^2)$ cost of orthonormalizing the columns of the matrix $Q$.

- $\tau_{A^*}(k+p)$ cost of forming matrix $B$ by multiplying $Q^*A$.

- $\mathcal{O}(n(k+p)^2)$ cost of the classical SVD.

- $\mathcal{O}(m(k+p)^2)$ cost of computing matrix $U$ as $Q\tilde{U}$.

Adding them together, gives the total cost of the general randomized SVD algorithm as

$$\tau_A(k+p) + \mathcal{O}(n(k+p)^2).$$

## 3.2  Power Iteration

Following the above analysis, it is easy to observe that the general randomized SVD method can perform badly if the singular values of the input matrix decay slowly. To solve this problem, we can combine Algorithm 2 with a few steps of power iteration [7]. Let $q$ be the number of steps of a power iteration, i.e., $q = 1$ or $q = 2$ in most cases.

Let $A = U\Sigma V^*$ be the singular value decomposition of matrix $A$. By running $q$ steps of power iteration, the singular values of matrix $A^{(q)}$ (defined below) are the diagonal entries on matrix

$\Sigma^{2q+1}$:

$$\begin{aligned} A^{(q)} &= (AA^*)^q A \\ &= ((U\Sigma V^*)(U\Sigma V^*)^*)^q A \\ &= (U\Sigma V^* V\Sigma U^*)^q A \\ &= (U\Sigma^2 U^*)^q A. \end{aligned} \tag{3.2}$$

Therefore, we see that the singular spectrum decays exponentially with $q$ steps of power iteration. Comparing $A^{(q)} = U\Sigma^{2q+1} V^*$ and the factorization of original matrix $A = U\Sigma V^*$, since SVD is unique, we see $A^{(q)}$ shares the same left and right singular vectors as $A$.

---
**Algorithm 3** Accuracy Enhanced Randomized SVD [13]
---

*Inputs:* An $m \times n$ matrix $A$, a target rank $k$, an oversampling parameter $p$ and an integer $q$ which is the number of steps in power iteration.

*Outputs:* Matrices U, V, and $\Sigma$ in an approximate rank-$(k+p)$ SVD of $A$. Where $U$ and $V$ are orthonormal, $\Sigma$ is diagonal.

**Stage 1:**
1: Generate an $n \times (k+p)$ Gaussian sketching matrix $\Omega$.
2: Form $Y = A\Omega$.
3: Apply $q$ steps of power iteration
4: **for** $j = 1 : q$ **do**
    $Z = A^* Y$
    $Y = AZ$
5: Construct a matrix Q whose columns form an orthonormal basis for the range of $Y$.

**Stage 2:**
6: Form $B = Q^* A$.
7: Compute an SVD of the small matrix: $B = \tilde{U}\Sigma V^*$.
8: Set $U = Q\tilde{U}$.

---

Let us now provide an analysis of Algorithm 3.

**Theorem 3** *Suppose that A is a real $m \times n$ matrix. Select an exponent q (e.g. $q = 1$ or 2), a target rank k and an oversampling parameter p, where $k + p \leq \min\{m, n\}$. Execute Algorithm 3 with a standard Gaussian test matrix $\Omega$ to obtain an $m \times (k+p)$ matrix Q with orthonormal columns.*

*Then*

$$\mathbb{E}\|A - QQ^*A\| \leq \left[1 + \frac{4\sqrt{(k+p)\min\{m,n\}}}{p-1}\right]^{1/(2q+1)} \sigma_{k+1},$$

$\sigma_{k+1}$ *is the (k+1)th largest singular value of A.*

Compare to Theorem 2, we can see that everything in the bracket is the same. Power iteration drives the leading constant to one exponentially fast as the exponent $q$ increases.

**Assumption 4** *[15] Let's make an assumption based on the Stage 1 of Algorithm 3. Let $\Omega_1 \in \mathbb{C}^{k \times (k+p)}$ by taking the partition of the matrix $V^*\Omega$*

$$V^*\Omega = \begin{bmatrix} V_k^*\Omega \\ V_\perp^*\Omega \end{bmatrix} = \begin{bmatrix} \Omega_1 \\ \Omega_2 \end{bmatrix},$$

*where $\Omega$ is a random Gaussian matrix from Stage 1 of Algorithm 3. We also assume that*

$$rank(\Omega_1) = k.$$

*The singular value gap at index k is inversely proportional to the singular value ratio*

$$\gamma_k = \|\Sigma_\perp\|_2 \|\Sigma_k^{-1}\|_2 = \frac{\sigma_{k+1}}{\sigma_k} < 1.$$

*Where $\Sigma_k$ is the Diagonal matrix that contains the first kth largest singular values of matrix A, and $\Sigma_\perp$ contains the rest of singular values of A.*

From this assumption, the first equation, $rank(\Omega_1 = k)$ guarantees that the starting guess $\Omega$ has a significant influence over the right singular vectors, the second inequality $\gamma_k < 1$ ensures that the $k$-dimensional subspace range ($U_k$) and also the right singular vectors formed $k$-dimensional subspace range ($V_k$) are both well-defined. In practice, we want $\gamma_k \ll 1$, so that there is a large singular value gap. With these assumptions, we can state the following theorem introduced in [15].

**Theorem 5** *[15] Let $U$, $V$ be obtained from Algorithm 3. Define the the canonical angles $\theta_j = \angle(U_k, U)$ and $v_j = \angle(V_k, V)$. Then under the Assumption 4, $\theta_j$ and $v_j$ satisfy*

$$\sin(\theta_j) \leq \frac{\gamma_j^{2q+1} \|\Omega_2 \Omega_1^\dagger\|_2}{\sqrt{1 + \gamma_j^{4q+2} \|\Omega_2 \Omega_1^\dagger\|_2^2}} \qquad \sin(v_j) \leq \frac{\gamma_j^{2q+2} \|\Omega_2 \Omega_1^\dagger\|_2}{\sqrt{1 + \gamma_j^{4q+4} \|\Omega_2 \Omega_1^\dagger\|_2^2}}$$

*for $j = 1, ..., k$*

## 3.3 Orthonormalization

Algorithm 3 now makes the randomized SVD method more compatible with the matrix that has singular spectrum decays slowly. But its accuracy can still be further improved due to round-off errors. While the power q increases, the columns of the sample matrix mentioned in algorithm 3,

$$Y = A^{(q)} \Omega \tag{3.3}$$

, tend to get closer and closer to the dominant left singular vector. Which causes all information that based on smaller singular values to get lots to round-off errors. [13]

And this can be solved by adding orthonormalization between each step of power iteration. This will make the algorithm more costly but at the same time more accurate.

**Algorithm 4** Accuracy Enhanced Randomized SVD with Orthonormalization[13]
___
*Inputs:* An $m \times n$ matrix $A$, a target rank $k$, an over-sampling parameter $p$ and an integer $q$ which is the number of steps in power iteration.

*Outputs:* Matrices $U$, $V$, and $\Sigma$ in an approximate rank-$(k+p)$ SVD of $A$. Where $U$ and $V$ are orthonormal, $\Sigma$ is diagonal.

**Stage 1:**

1: Generate an $n \times (k+p)$ Gaussian test matrix $\Omega$.

2: Form $Y = A\Omega$.

3: Construct a matrix Q whose columns form an orthonormal basis for the range of Y.

4: Apply power iterations:

5: **for** $j = 1 : q$ **do**

   $W = orth(A^*Q)$ ;

   $Q = orth(AW)$ ;

**Stage 2:**

6: Form $B = Q^*A$.

7: Compute an SVD of the small matrix: $B = \tilde{U}\Sigma V^*$.

8: Set $U = Q\tilde{U}$.
___

Based on Algorithm 4, we can introduce the following error analysis.

**Theorem 6** *[13, Theorem 2] Suppose that A is an m × n matrix. Select an exponent q (e.g. q = 1 or 2), a target rank and an oversampling parameter p, where $k + p \leq \min\{m, n\}$. Draw a Gaussian matrix $\Omega$ of size $n \times (k+p)$, define $Y = A^{(q)}\Omega$ and Q to be the $m \times (k+p)$ orthonormal matrix resulting from orthonormalizing the columns of Y. Then*

$$\mathbb{E}\|A - QQ^*A\| \leq \left[\left(1 + \sqrt{\frac{k}{p-1}}\right)\sigma_{k+1}^{2q+1} + \frac{e\sqrt{k+p}}{p}\left(\sum_{j=k+1}^{min(m,n)} \sigma_j^{2(2q+1)}\right)^{1/2}\right]^{1/(2q+1)}, \qquad (3.4)$$

*where $\sigma_{k+1}$ is the (k+1)th largest singular value of A.*

We can simplify this result by considering the worst case where there is no decay in the singular values after the $k$th term, i.e., $\sigma_{k+1} = \sigma_{k+2} = ... = \sigma_{min(m,n)}$. Then,

$$\mathbb{E}\|A - QQ^*A\| \leq \left[\left(1 + \sqrt{\frac{k}{p-1}}\right) + \frac{e\sqrt{k+p}}{p} \cdot \sqrt{\min\{m,n\} - k}\right]^{1/(2q+1)} \sigma_{k+1}. \quad (3.5)$$

Now as the exponent $q$ increases, the power iteration scheme drives the factor in front of $\sigma_{k+1}$ to one exponentially fast.

## 3.4 Single Pass

Randomization can also improve efficiency from the storage side.

All randomized SVD algorithms we have introduced so far require the access to the large input matrix $A$ twice, for both Stage 1 and Stage 2. It is possible for us to modify the algorithm such that each entry of matrix $A$ is accessed only once, so a significantly amount of storage can be reduced.

### 3.4.1 Hermitian Matrices

In the case of the Hermitian input matrix, a single-pass algorithm can be directly developed from Algorithm 2. For the single pass randomized SVD method with Hermitian matrix Stage 1 is the same as in Algorithm 1. We first draw an $n \times (k+p)$ Gaussian random matrix $\Omega$ and form the matrix $Y = A\Omega$. Then we construct a matrix $Q$ whose columns form an orthonormal basis for the range of $Y$. This yields

$$A \approx QQ^*A. \quad (3.6)$$

Given that $A$ is a Hermitian matrix, we have that $A = A^*$, therefore

$$A \approx AQQ^*. \quad (3.7)$$

Replacing *A* in (3.7) with (3.6) gives

$$A \approx QQ^*AQQ^*. \tag{3.8}$$

Setting

$$C := Q^*AQ, \tag{3.9}$$

allows us to apply the eigenvalue decomposition of matrix *C* to find $\hat{U}$ and *D* such that $C = \hat{U}D\hat{U}^*$. With $U := Q\hat{U}$ we get the following result

$$A \approx QCQ^* = Q\hat{U}D\hat{U}^*Q^* = UDU^*. \tag{3.10}$$

Since our goal is to build a single pass algorithm, we cannot compute matrix *C* explicitly. Instead, we will multiply both side of equation (3.9) by $Q^*\Omega$ on the right to obtain

$$C(Q^*\Omega) = Q^*AQ(Q^*\Omega). \tag{3.11}$$

Now by equation (3.7),

$$C(Q^*\Omega) = Q^*(AQQ^*)\Omega \approx Q^*A\Omega = Q^*Y. \tag{3.12}$$

If we ignore the approximation error, we can compute matrix *C* as the solution of the linear system

$$C(Q^*\Omega) = Q^*Y. \tag{3.13}$$

This leads to our single-pass randomized SVD method for a Hermitian matrix given in Algorithm 5.

Since computing matrix *C* resulting additional approximation error, this single-pass method is less accurate compared to previously discussed two-pass methods.[13]

**Algorithm 5** Single-Pass Randomized SVD for a Hermitian Matrix [13]

    *Inputs:* An $n \times n$ Hermitian matrix $A$, a target rank $k$, an oversampling parameter $p$.
    *Outputs:* Matrices $U$, and $D$ in an approximate rank-$k$ EVD of $A$. Where $U$ is orthonormal, $\Sigma$ is diagonal.
    **Stage 1:**
1: Generate an $n \times (k+p)$ Gaussian test matrix $\Omega$.
2: Form $Y = A\Omega$.
3: Let $Q$ denote the orthonormal matrix formed by $Y$.
    **Stage 2:**
4: Let $C$ denote the $k \times k$ solution of $C(Q^*\Omega) = (Q^*Y)$
5: Compute an eigenvalue decomposition of the small matrix $C = \widehat{U}D\widehat{U}^*$.
6: Set $U = Q\widehat{U}$.

### 3.4.2 General Matrices

Algorithm 5 works only for Hermitian input matrix limiting its use. Let us now consider the single-pass idea in the case of general input matrices. First, in Stage 1, we need to apply randomized sampling simultaneously to both the row and the column space of matrix $A$ [13]. We start with generating two random Gaussian matrices $\Omega_c \in \mathbb{C}^{n \times (k+p)}$ and $\Omega_r \in \mathbb{C}^{m \times (k+p)}$, and forming two sketches $Y_c = A\Omega_c$ and $Y_r = A^*\Omega_r$. Then we construct two matrices $Q_c = orth(Y_c)$ and $Q_r = orth(Y_r)$. For Stage 2, we execute a projection step to obtain a smaller matrix

$$C = Q_c^* A Q_r. \tag{3.14}$$

To develop a single-pass algorithm, we start with left multiplying (3.14) by $\Omega_r^* Q_c$, i.e.,

$$\Omega_r^* Q_c C = \Omega_r^* Q_c Q_c^* A Q_r \approx \Omega_r^* A Q_r = Y_r^* Q_r. \tag{3.15}$$

Then, we right multiply (3.14) by $Q_r^* \Omega_c$, i.e.,

$$C Q_r^* \Omega_c = Q_c^* A Q_r Q_r^* \Omega_c \approx Q_c^* A \Omega_c = Q_c^* Y_c. \tag{3.16}$$

This allows us to obtain $C$ as the least-square solution of the two equations

$$(\Omega_r^* Q_c)C = Y_r^* Q_r \quad \text{and} \quad C(Q_r^* \Omega_c) = Q_c^* Y_c.$$

All discussed steps form our single-pass randomized SVD method for general matrices presented in Algorithm 6.

---
**Algorithm 6** Single-Pass Randomized SVD for a General Matrix [13]
---
*Inputs:* An $m \times n$ matrix $A$, a target rank $k$, an oversampling parameter $p$.
*Outputs:* Orthonormal matrices U, V, and diagonal $\Sigma$ in an approximate rank-$k$ SVD of $A$.
**Stage 1:**
1: Generate two Gaussian matrices $\Omega_c$ and $\Omega_r$ of size $n \times (k+p)$ and $m \times (k+p)$, respectively.
2: Form $Y_c = A\Omega_c$, $Y_r = A\Omega_r$.
3: Construct orthonormal matrices $Q_c$ and $Q_r$ consisting of the $k$ dominant left singular vectors of $Y_c$ and $Y_r$.
**Stage 2:**
4: Compute C as the solution of the joint system of equations formed by $(\Omega_r^* Q_c)C = Y_r^* Q_r$ and $C(Q_r^* \Omega_c) = Q_c^* Y_c$
5: Compute an SVD of the small matrix: $C = \widehat{U}\Sigma\widehat{V}^*$.
6: Set $U = Q_c\widehat{U}$ and $V = Q_r\widehat{V}$.

---

# Chapter 4

# Numerical Experiments

## 4.1 Illustration of various randomized SVD algorithms

First, we are going to use image reconstruction example to show the performance of all methods discussed in Chapter 3.

Figure 4.1 is the original image that we are going to use in this experiment. This figure is converted to a 804 × 1092 matrix of rank 804.

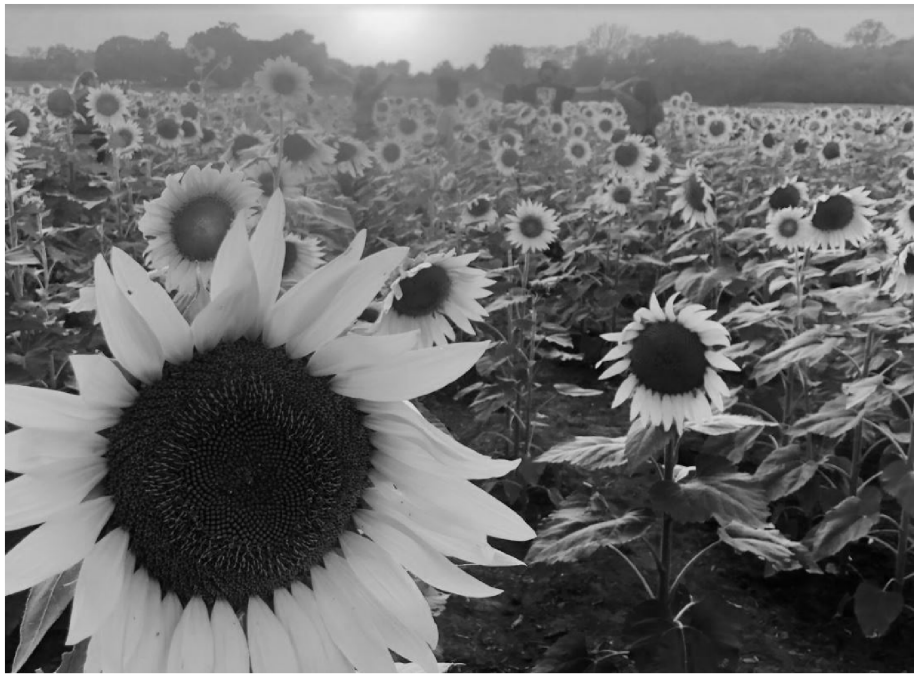

Figure 4.1: Original Sunflower image of rank 804.

### 4.1.1 Algorithm 2 - General Randomized SVD Method

First, Figure 4.2 presents rank-$k$ approximations of the Sunflower image with fixed $p = 0$ and $k = 10, 50, 100, 400, 800$ obtained using Algorithm 2. Figure 4.3 present the rank-$k$ approximation obtained with fixed $p = 10$, $k = 10, 50, 100, 400, 800$.

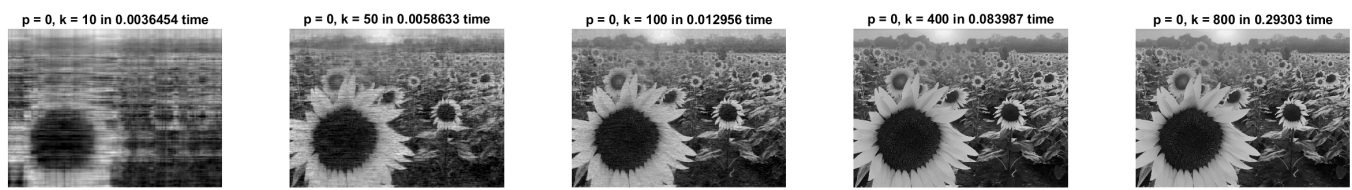

Figure 4.2: Illustration of low-rank approximation of Sunflower image using Algorithm 2 with no oversampling $p = 0$ and $k = 10, 50, 100, 400, 800$ (left to right).

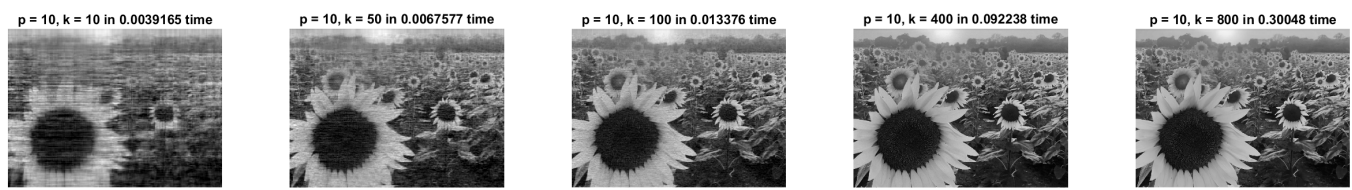

Figure 4.3: Illustration of low-rank approximation of Sunflower image using Algorithm 2 with oversampling $p = 10$ and $k = 10, 50, 100, 400, 800$ (left to right).

Comparing these two figures horizontally, we see that for the same value of parameter $p$ the larger target rank $k$ the better the quality of the approximation, i.e., the restored images are more clear. Then comparing vertically, we observe that with the same target rank $k$, the larger the oversampling parameter $p$ is, the more clear the restored picture.

### 4.1.2 Algorithm 3 - Accuracy Enhanced Randomized SVD Method

Similar experiment was performed with Algorithm 3. We first fix the number of power iterations $q$ involved in Step 1 and let $k = 10, 50, 100, 400, 800$. The corresponding results are presented in Figure 4.4 for $q = 1$, and Figure 4.5 for $q = 2$.

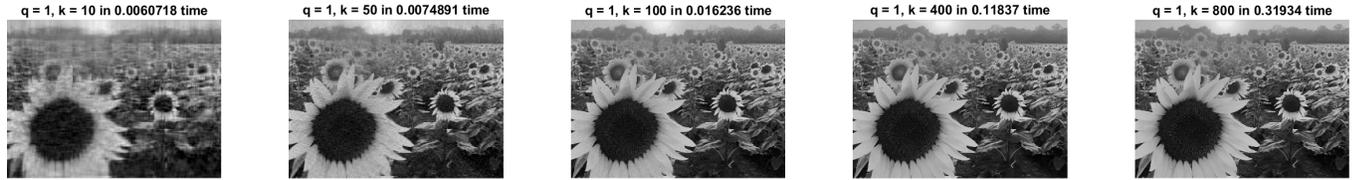

Figure 4.4: Illustration of low-rank approximation of Sunflower image using Algorithm 3 with $q = 1$ steps of power iteration and $k = 10, 50, 100, 400, 800$ (left to right).

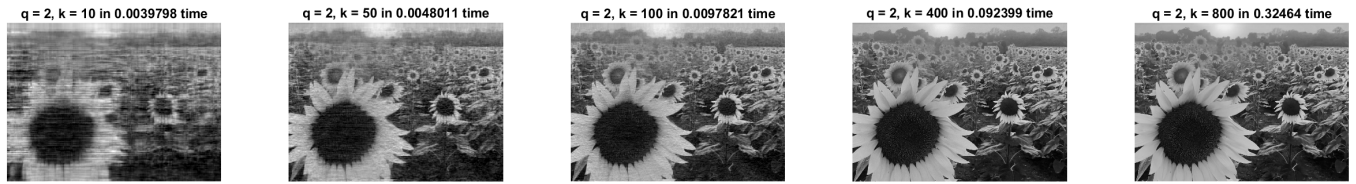

Figure 4.5: Illustration of low-rank approximation of Sunflower image using Algorithm 3 with $q = 2$ steps of power iteration and $k = 10, 50, 100, 400, 800$ (left to right).

Comparing these two figures horizontally, we see that with the same value of parameter $q$, when the target rank $k$ gets larger the restored picture is more clear.

### 4.1.3 Algorithm 4 - Accuracy Enhanced Randomized SVD with Orthonormalization

Let us now repeat the above experiment with Algorithm 4. Figure 4.6 for $q = 1$, and Figure 4.7 present the results for $q = 1$ and $q = 2$, respectively.

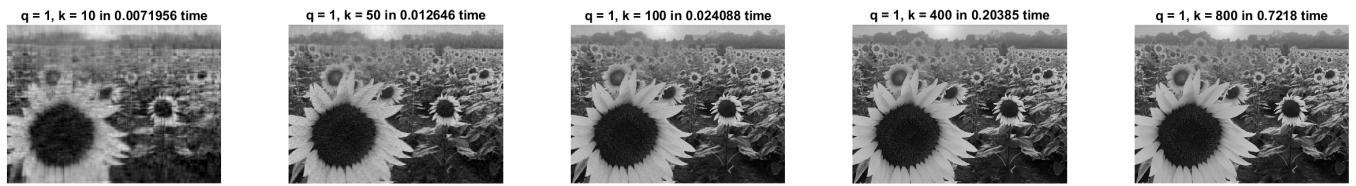

Figure 4.6: Illustration of low-rank approximation of Sunflower image using Algorithm 4 with $q = 1$ steps of power iteration and $k = 10, 50, 100, 400, 800$ (left to right).

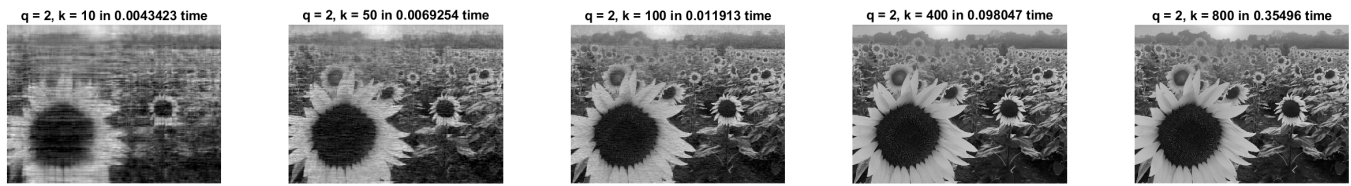

Figure 4.7: Illustration of low-rank approximation of Sunflower image using Algorithm 4 with $q = 2$ steps of power iteration and $k = 10, 50, 100, 400, 800$ (left to right).

Comparing these two figures horizontally, we see that with the same value of parameter $q$, when the target rank $k$ gets larger, the restored pictures are more clear.

## 4.2 Illustration of Error Bounds

Our next experiment will illustrate the theoretical error bounds presented in Theorem 2, 3 and 4 derived for the Stage 1 of Algorithm 2, 3 and 4, respectively, and compare them with the computed error $\|A - QQ^*A\|_2$. Let us consider a set of test matrices $A$ defined as in [15, Section 6.1]. We take a set of sparse matrices $A \in \mathbb{R}^{3000 \times 300}$ with prescribed decay of singular values, location and size/value of the singular values gap. Given a sparse vectors $x_j \in \mathbb{R}^{3000}$ and $y_j \in \mathbb{R}^{300}$ with nonnegative entries, the location $r$ of the gap between the singular values of a test matrix $A$ and its size determined by the parameter $gap$ we define $A$ according to the following formula

$$A = \sum_{j=1}^{r} \frac{\text{gap}}{j} x_j y_j^T + \sum_{j=r+1}^{300} \frac{1}{j} x_j y_j^T. \tag{4.1}$$

Notice that the singular values of such constructed matrix $A$ decay like $\frac{1}{j}$.

First, let us consider three test matrices $A$ each with a different size of the $gap$ parameter, i.e., small, medium and large gap corresponding to the value $gap = 1, 2$ or $10$, respectively, between the 15th and 16th largest singular values, i.e., $r = 15$. Figure 4.8 shows the exact singular values of each of the three constructed test matrices.

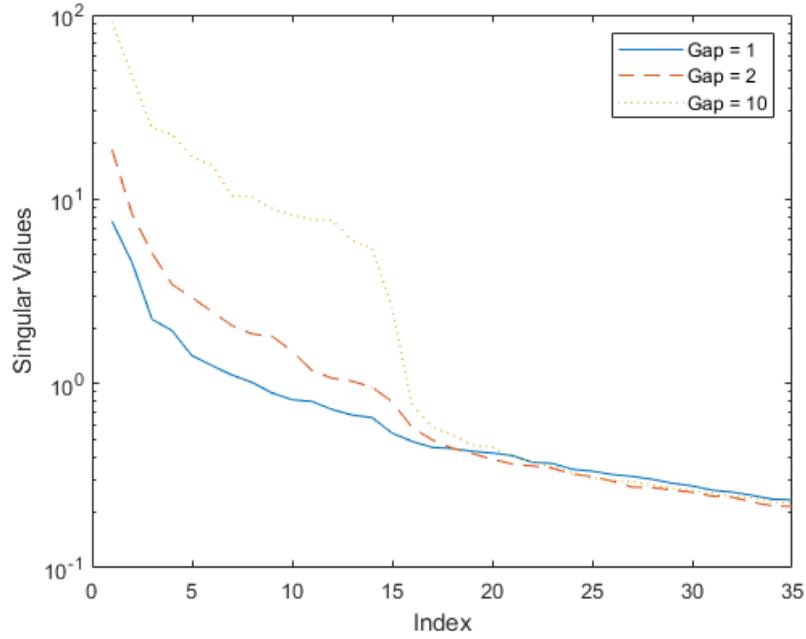

Figure 4.8: Singular values of the three test matrices defined in (4.1) with parameter $gap = 1, 2$ and 10.

### 4.2.1 General Randomized SVD - Fixed Rank Problem

Let us now fix the size of the $gap$ to be 10 and take the corresponding test matrix $A$ from (4.1) to be the input of the Algorithm 1. Since Theorem 2 provides a theoretical upper bound of the expectation value of the $2-norm$ of the difference between the original input matrix $A$ and its low-rank approximation, i.e., $\|A-QQ^*A\|_2$, in the following paragraphs we will compare this theoretical bound (the estimated error) with the size of the actual error $\|A-QQ^*A\|_2$ (computed error).

**Error analysis with fixed $p$ and different ranks $k$**  For this experiment, we fix the value of the oversampling parameter $p = 5$ choose different target ranks $k$ to be $5, 10, 15, 20, 25, 30$. Figure 4.9 illustrates the computed errors and estimated errors for Algorithm 1 with different values of $k$.

**Error analysis with fixed rank $k$ and different values of parameter $p$**  Here, we fix the value of the target rank $k = 20$ and choose different values of $p = 5, 10, 15, 20, 25$. Figure 4.10 illustrates the computed errors and estimated errors for Algorithm 1 with different values of parameter $p = 5, 10, 15, 20, 25$.

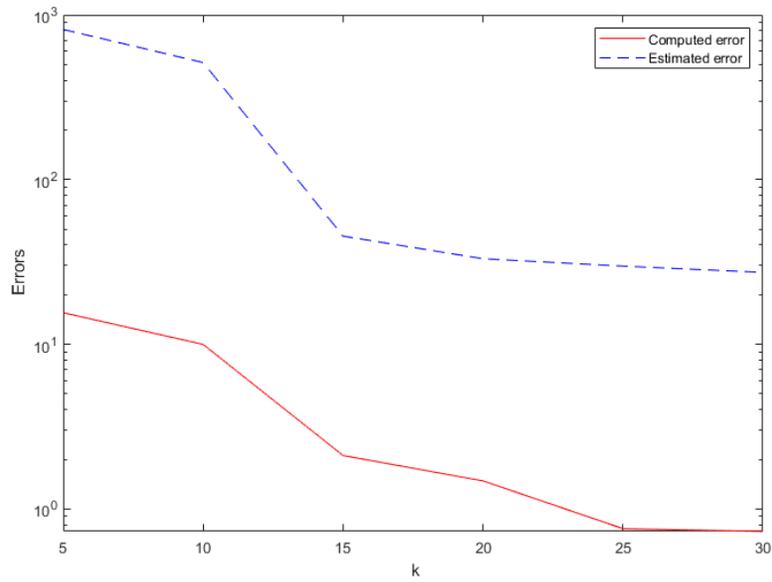

Figure 4.9: Illustration of error bounds for Algorithm 1 obtained for test matrix $A$ with fixed $p = 5$, and different values of rank $k$.

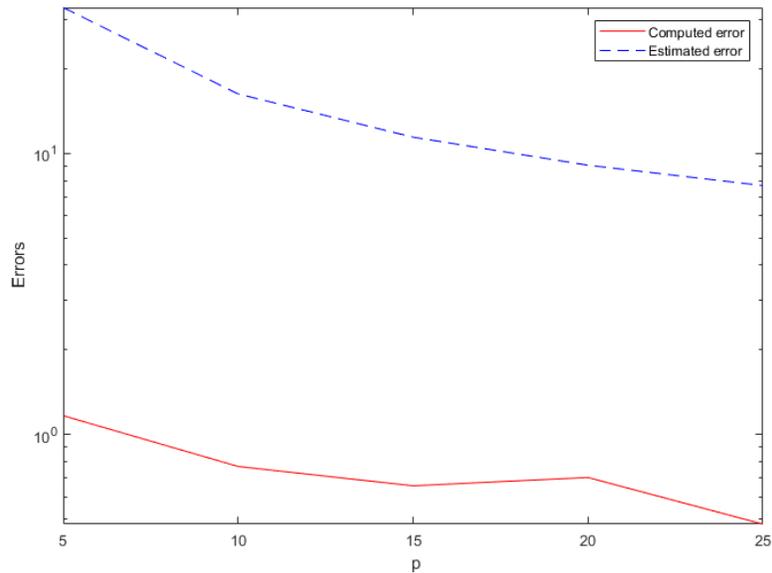

Figure 4.10: Illustration of error bounds for Algorithm 1 obtained for test matrix $A$ with fixed rank $k = 20$ and different values of parameter $p$.

**Average of computed errors over 100 runs** With fixed target rank $k = 20$ and oversampling parameter $p = 5$, Figure 4.11 presents an average computed and estimated errors over 100 runs of Algorithm 1.

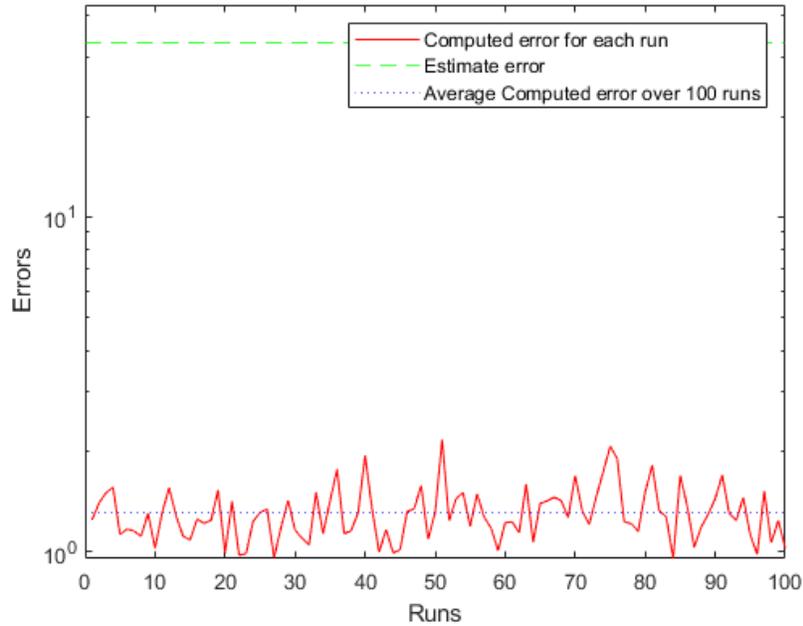

Figure 4.11: Illustration of average errors obtained for the test matrix $A$ with fixed rank $k = 20$ and oversampling parameter $p = 5$ over 100 runs of Algorithm 1.

### 4.2.2 Accuracy Enhanced Randomized SVD - Power Iteration

We use the same test matrix $A$ and the same experiment parameters as above for the accuracy enhanced randomized SVD method in Algorithm 3. Here, Stage 1 of the algorithm incorporates few steps, e.g., $q = 1$ or $2$ of power iteration method.

**Error analysis with fixed parameters $p$ and $q$, and different values of $k$** For this experiment, we fix the oversampling parameter $p = 5$ and the number of the power iteration steps $q = 1$, and vary the value of target rank parameter $k = 5, 10, 15, 20, 25, 30$. Figure 4.12 presents the computed error and estimated error of Algorithm 3 with different values of $k$.

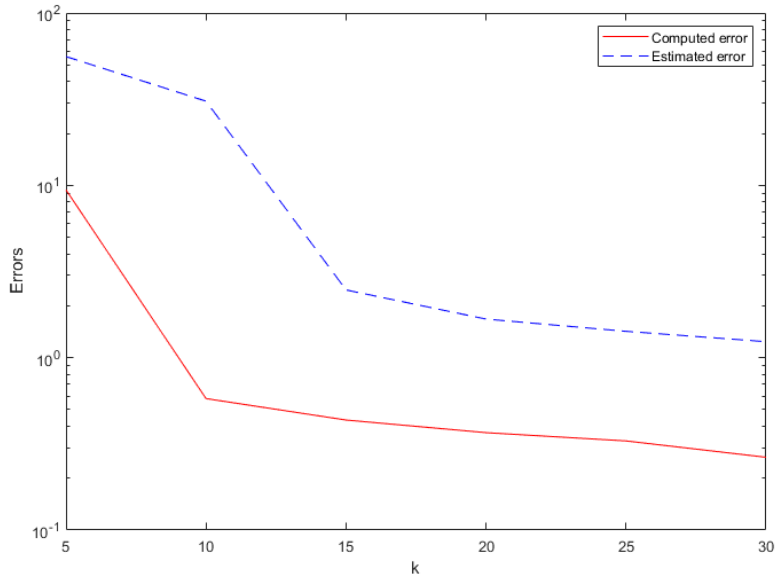

Figure 4.12: Illustration of error bounds for Algorithm 3 obtained for test matrix $A$ with fixed values of $p = 5$ and $q = 1$, and different values of $k$.

**Error analysis with fixed values of $k$ and $q$, and different values of parameter $p$**   For this experiment, we fix a target rank $k = 20$ and the number of power iteration steps $q = 1$ and choose different values of oversampling parameter $p = 5, 10, 15, 20, 25$. Figure 4.14 presents the computed error and estimated error of Algorithm 3 with different values of $p$.

**Average of computed error over $100$ runs**   With fixed rank $k = 20$, oversampling parameter $p = 5$ and number of power iteration steps $q = 1$, Figure 4.14 presents an average computed and estimated errors over 100 runs of Algorithm 3. Note that with the fixed values of parameters $k$, $p$ and $q$, the estimated error from Theorem 3 is constant.

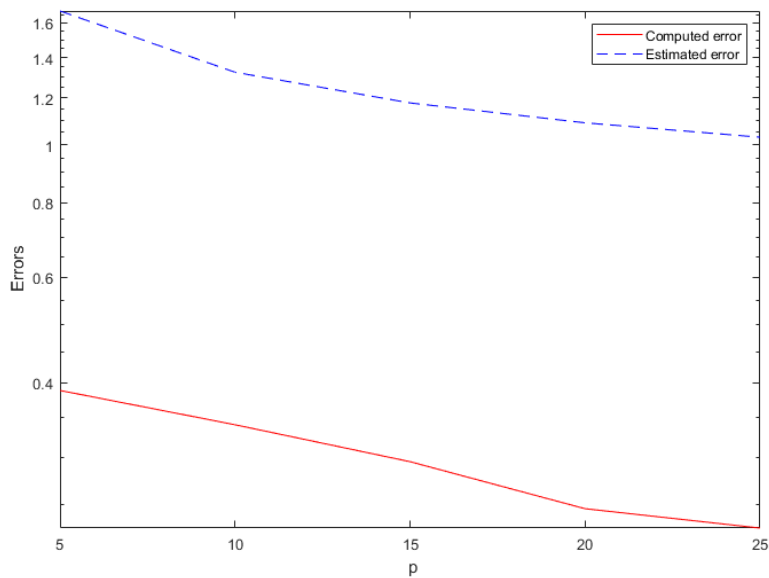

Figure 4.13: Illustration of error bounds for Algorithm 3 obtained for test matrix $A$ with fixed values of $k = 20$ and $q = 1$, and different values of $p$.

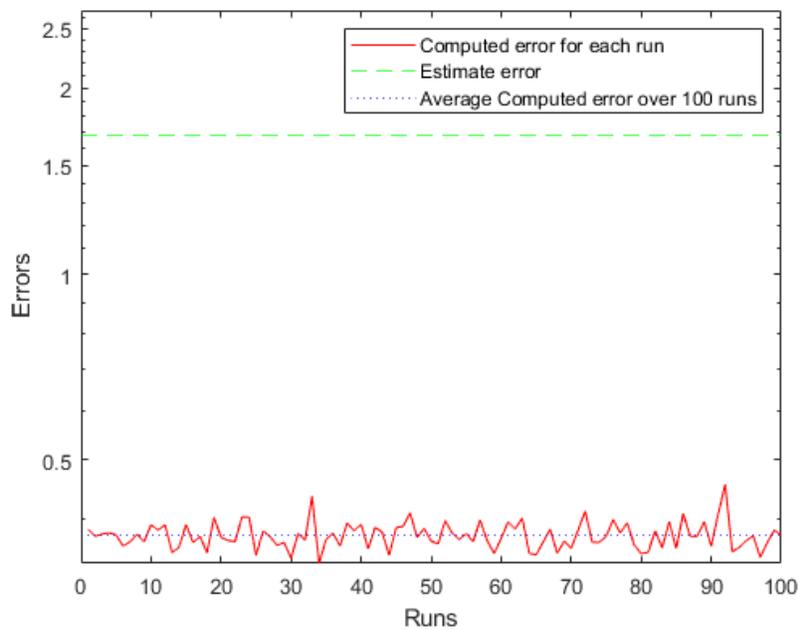

Figure 4.14: Illustration of average error bounds obtained for test matrix $A$ with $k = 20$, $p = 5$, and $q = 1$ over 100 runs of Algorithm 3.

### 4.2.3 Accuracy enhanced randomized SVD with Orthonormalization - Power Iteration

Similarly, for (Alg 4) the Accuracy enhanced randomized SVD with Orthonormalization, the following figures are generated.

**Error analysis with fixed $p$, $q$ and varies $k$**  For this experiment, we fix the oversampling parameter $p$ to be 5, $q$ to be 1 and choose $k$ to be 5, 10, 15, 20, 25, 30, respectively. And compute the computed error and estimated error for different $k$.

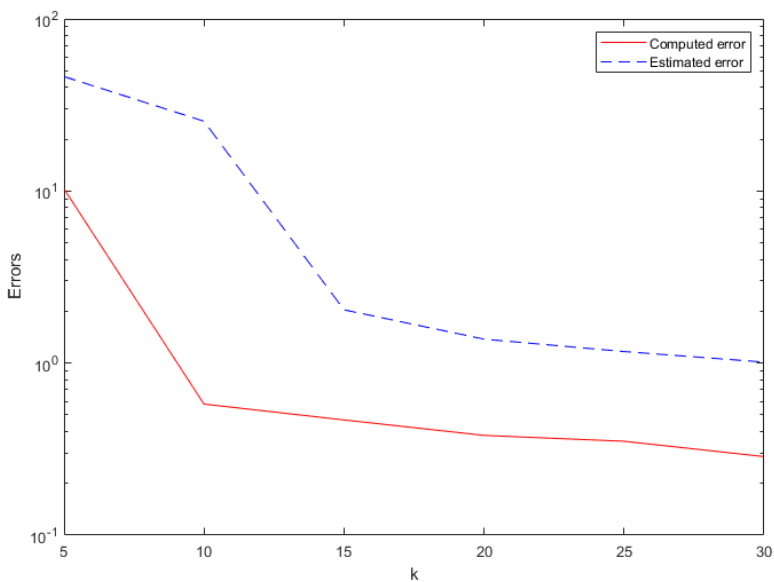

Figure 4.15: Illustration of error bounds of the test matrix $A$ with fixed $p = 5$, $q = 1$, and $k = 5, 10, 15, 20, 25, 30$ (left to right).

**Error analysis with fixed $k$, $q$ and varies $q$**  For this experiment, we fix the oversampling parameter $k$ to be 20, $q$ to be 1 and choose $p$ to be 5, 10, 15, 20, 25, respectively. And compute the computed error and estimated for different $p$.

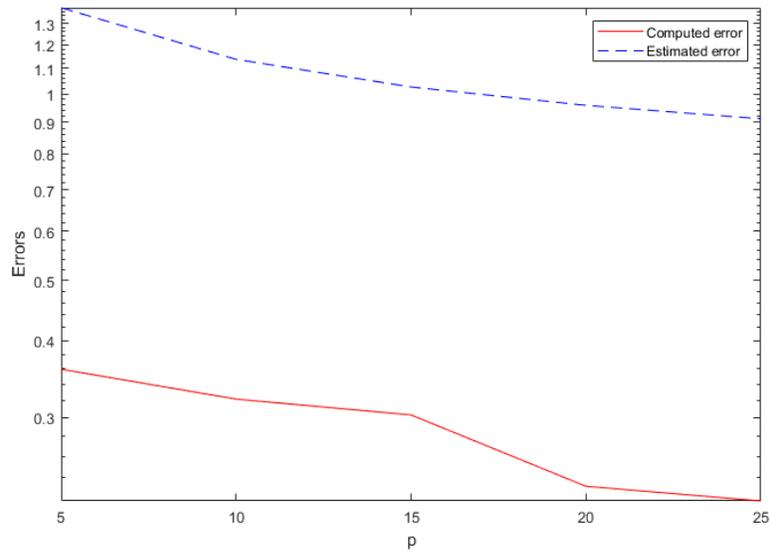

Figure 4.16: Illustration of error bounds of the test matrix $A$ with fixed $k = 20$, $q = 1$, and $p = 5, 10, 15, 20, 25$ (left to right).

**Average of computed error over** $100$ **runs** With fixed $k = 20$, $p = 5$, $q = 1$, this experiment run the stage 1 of algorithm 4 for 100 times and generates the average computed error. Note that with $k$, $p$ and $q$ fixed, the estimated error from theorem 6 is also fixed.

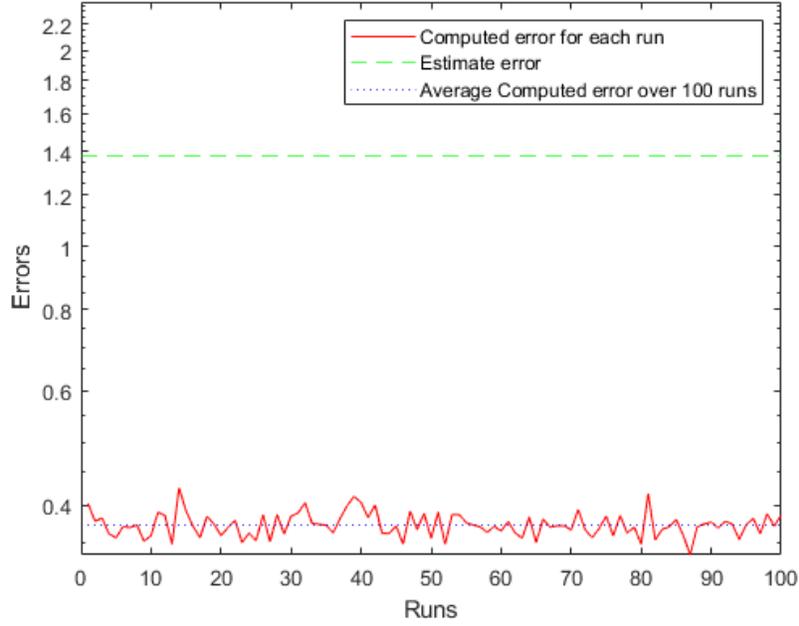

Figure 4.17: Illustration of average error bounds of the test matrix $A$ with fixed $k = 20$, $p = 5$, and $q = 1$ over 100 runs.

## 4.3  Illustration of Canonical angles

For this experiment we take the same test matrices $A$ as above with the $gap = 2$ between the 15th and 16th largest singular values. We consider the target rank $k = 25$ and the oversampling parameter $p = 5$.

Let us first recall the singular value decomposition of matrix $A$. Given a target rank $k$ we write the singular value decomposition of matrix $A$ as follows:

$$A = [U_k \quad U_\perp] \begin{bmatrix} \Sigma_k & \\ & \Sigma_\perp \end{bmatrix} \begin{bmatrix} V_k^* \\ V_\perp^* \end{bmatrix}. \tag{4.2}$$

Here, $\Sigma_k \in \mathbb{C}^{k \times k}$ and $\Sigma_\perp \in \mathbb{C}^{(m-k) \times (n-k)}$ are diagonal matrices; the columns of $U_k$ and $U_\perp$ are the corresponding left singular vectors, and the columns of $V_k$ and $V_\perp$ the corresponding right singular vectors of matrix $A$. We also denote by $A_k = U_k \Sigma_k V_k^*$ the best rank-$k$ approximation of matrix $A$.

Our goal is to determine how well range($\widehat{U}$) approximates the range($U_k$) and we will achieve this my measuring the canonical angles between these two subspace, i.e., calculating the

$$\sin \angle (U_k, \widehat{U}). \tag{4.3}$$

We will consider $\widehat{U}$ resulting from Algorithm 2, 3 and 4.

**Canonical angles computed by Algorithms 2, 3 and 4 with $q = 0$** Figure 4.18 shows the computed canonical angles for Algorithm 2, 3 and 4 applied to the test matrix $A$ and $q = 0$ steps of power iteration. We see that since $q = 0$ to Algorithm 3 and 4 did not perform any steps of power iteration in Stage 1 which makes them easily comparable to Algorithm 2 which first solves the general fixed rank problem, see Algorithm 1. This observation is validated with the results presented on Figure 4.18.

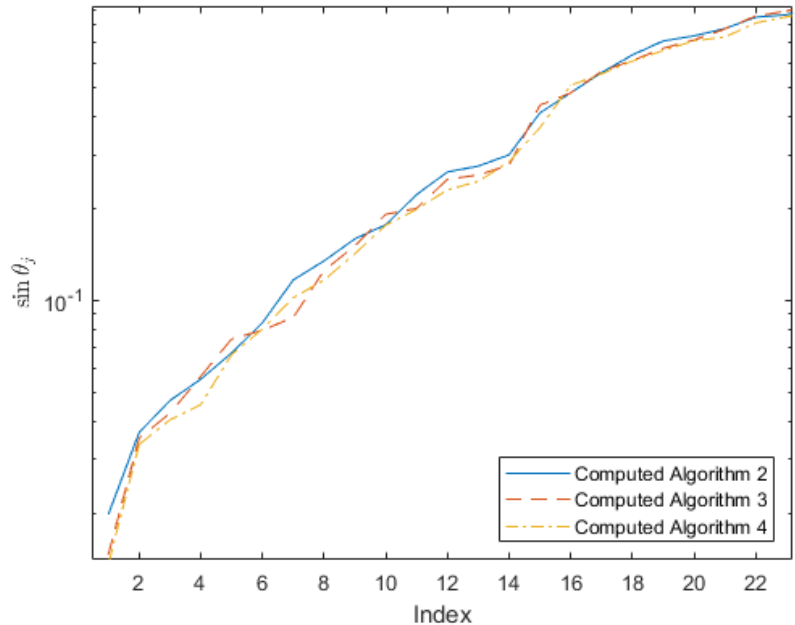

Figure 4.18: Illustration of the canonical angles determined for matrix $A$ through Algorithm 2, 3 and 4 with $q = 0$.

**Computed and estimated canonical angles for Algorithm 3**  We plot the canonical angles in solid lines, and the corresponding bounds from Theorem 5 in dashed lines. The results are based on Algorithm 3 with input matrix $A$ and the steps of a power iteration to be $q = 0, 1, 2$:

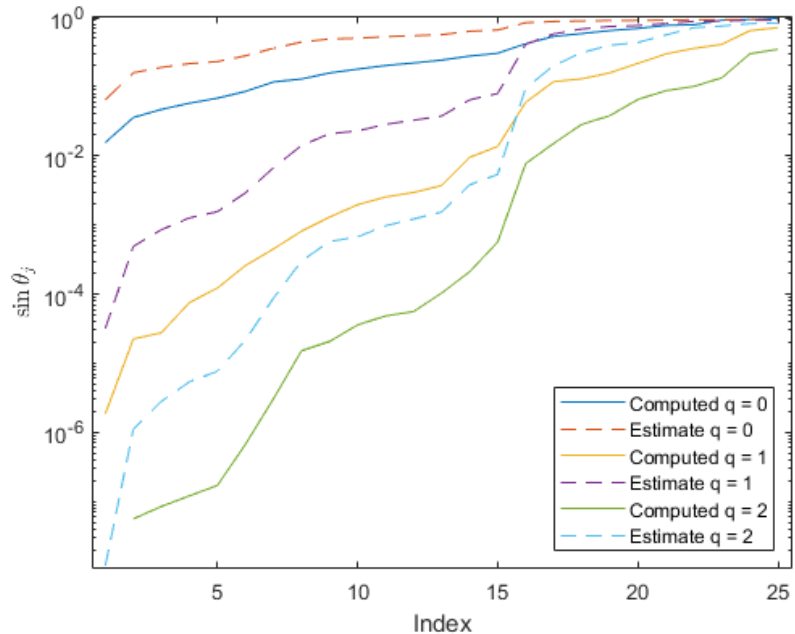

Figure 4.19: Illustration of canonical angles of the test matrix $A$ for Alg 3 with $q = 0, 1, 2$. The solid lines correspond to the computed values, the dashed lines correspond to bounds obtained using Theorem 5

From the figure 4.19, we observe that when $q$ increases, the size of the canonical angle becomes smaller. So the 'sketch' is becoming more accurate. Then within the same value $q$, as the index increases, the canonical angle increases. As we setup the matrix $A$ with a size 2 $gap$ between 15th and 16th singular values, all the canonical angles below index 15 are captured accurately.

**Computed and estimated canonical angles for Alg 4**  Now we apply the same experiment above with algorithm 4. We still use the same estimate bounds stated in Theorem 5, since by our construction of Algorithm 4, we should expect a similar error analysis but more accurate when $q$ is large, since algorithm 4 is designed to limit the round-off errors in algorithm 3.

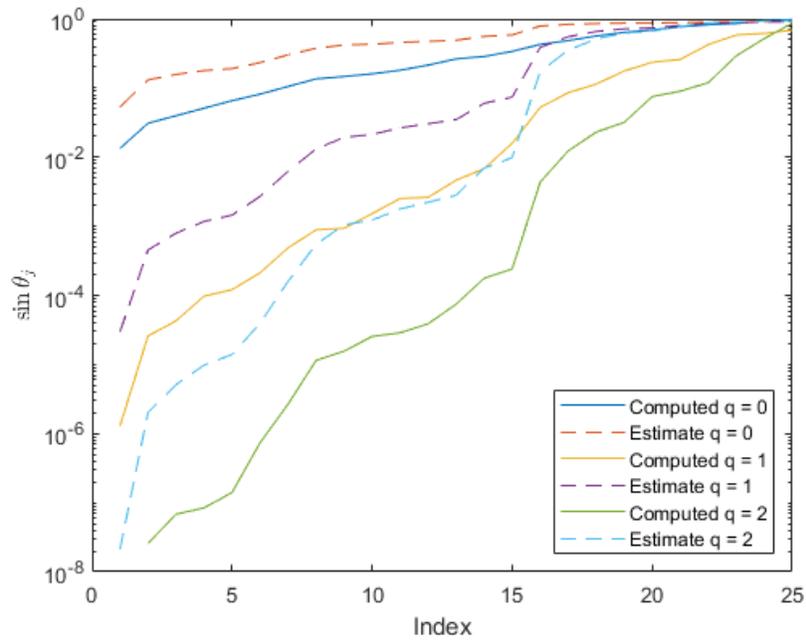

Figure 4.20: Illustration of canonical angles of the test matrix $A$ for Alg4 with $q = 0, 1, 2$. The solid lines correspond to the computed values, the dashed lines correspond to bounds obtained using Theorem 5

We are getting the same (which is what we expected) observation as from last experiment. If we compare this figure 4.20 with last figure 4.19, we can also see that the difference between the computed error and the estimate error is later from this figure when $q = 2$(large). Which confirms the expectation that algorithm 4 reduces the round off error that algorithm 3 has when $q$ is large.

# Chapter 5

# Conclusions

In this report, we present different versions of randomized algorithms for calculating the singular value decomposition and the corresponding theoretical error bounds. We first introduce the general randomized SVD and then consider incorporating few steps of power iteration within the Stage 1 to make the algorithms efficient when the singular spectrum of the input matrix decays slowly. Orthonormalization can be involved in each power iteration step to reduce the round-off error and further improve the accuracy. On the other hand, we can save the storage costs making the randomized SVD algorithm only access each entry of input matrix $A$ once. We first introduce this single-pass method for the Hermitian input matrix $A$ and then extend it to work for general matrices.

In the first (image restoration) experiment, we illustrate the performance of various randomized SVD methods and observe that when the value of target rank $k$ and/or oversampling parameter $p$ increases, the reconstructed images are more and more accurate.

The second experiment compares the theoretical error bounds with the actually computed errors. We can see that our computed errors are always below the theoretical bounds, and both estimated and computed errors decrease when the target rank $k$ and/or oversampling parameter $p$ increase.

Finally, the third experiment illustrates the quality of approximations obtained with randomized SVD in terms of canonical angles. We observe that when the number of power iterations steps $q$ is larger, the canonical angles get smaller, confirming that the corresponding random 'sketch' provides good approximation of the original matrix $A$.

Future work may include more experiments with cross-comparisons between the error per-formance of Algorithm 3 and Algorithm 4 when the number of power iteration steps $q$ gets large. Since Algorithm 4 is designed to reduce the round-off error with respect to Algorithm 3, we should expect to see significantly lower error for Algorithm 4 for large values of $q$.

# Appendix A

# Randomized NLA MATLAB Toolbox

**Code Link:** As a part of this project, we have created a simple randomized numerical linear algebra MATLAB Toolbox which can be found on GitHub https://github.com/LeeeeeLy/Randomized-NLA-Matlab-Toolbox

In the following we will briefly describe the structure of our toolbox.

**Main Functions:**

- **FixedRank.m** - This is an implementation of Proto-Algorithm: Solving the Fixed-Rank Problem from paper [11].

  This function takes an input matrix $A$, a target rank $k$, and a oversampling parameter $p$ and output a $m \times (k+p)$ matrix $Q$ whose columns are orthonormal and whose range approximates the range of $A$.

- **RandSVD.m** - This is implementation of Prototype for Randomized SVD from paper [11].

  This function takes an input matrix $A$, a target rank $k$, and an exponent $q$ ($q$ = 0,1,2) and approximate a rank-$2k$ factorization $A \approx U\Sigma V^*$, where $U$ and $V$ are orthonormal, and $\Sigma$ is nonnegative and diagonal.

- **randRF.m** - This is implementation of Algorithm 4.1 from paper [11].

  This function takes an $m \times n$ matrix $A$ and an integer $l$, it computes an $m \times l$ orthonormal matrix $Q$ whose range approximates the range of $A$.

- **Adaptive_randRF.m** - This is implementation of Algorithm 4.2 from paper [11].

  This function takes an $m \times n$ matrix $A$, a tolerance $\epsilon$, and an integer $r$ (e.g., $r = 10$), it computes an orthonormal matrix $Q$ such that $\|(I - QQ^*)A\| \leq \epsilon$, with probability at least $1 \min\{m, n\} 10^{-r}$.

- **randPI.m** - This is implementation of Algorithm 4.3 from paper [11].

  This function takes an $m \times n$ matrix $A$ and integers $l$ and $q$, it computes an $m \times l$ orthonormal matrix $Q$ whose range approximates the range of $A$.

- **randSI.m** - This is implementation of Algorithm 4.4 from paper [11].

  This function takes an $m \times n$ matrix $A$ and integers $l$ and $q$, it computes an $m \times l$ orthonormal matrix $Q$ whose range approximates the range of $A$.

- **FastRandRF.m** - This is implementation of Algorithm 4.5 from paper [11].

  This function takes an $m \times n$ matrix $A$ and an integer $l$, it computes an $m \times l$ orthonormal matrix $Q$ whose range approximates the range of $A$. This function involved with sub sampled random Fourier transform (SRFT), see function **SRFT.m** in the Sub-Function session.

- **DirectEigvalueDecopo.m** - This is implementation of Algorithm 5.3 from paper [11].

  This function takes an Hermitian matrix $A$ and a basis $Q$ that can be generated by **FixedRank.m**, this computes an approximate eigenvalue decomposition $A \approx U \Lambda U$, where $U$ is orthonormal, and $\Lambda$ is a real diagonal matrix.

  To generate an input hermitian matrix, please visit and download the [random hermitian matrix generator function](#)[12].

- **EigvalueDecopoRow.m** - This is implementation of Algorithm 5.4 from paper [11].

  This function takes an Hermitian matrix $A$ and a basis $Q$ that can be generated by FixedRank.m, this computes an approximate eigenvalue decomposition $A \approx U \Lambda U^*$, where $U$ is orthonormal, and Λ is a real diagonal matrix. This function is faster than **DirectEigvalueDecopo.m** but less accurate.

To generate an input hermitian matrix, please visit and download the [random hermitian matrix generator function](#)[12].

- **EigvalueDecopoNystrom.m** - This is implementation of Algorithm 5.5 from paper [11].

  This function takes a positive semidefinite matrix $A$ and a basis $Q$ that can be generated by FixedRank.m, this computes an approximate eigenvalue decomposition $A \approx U\Lambda U^*$, where $U$ is orthonormal, and Λ is nonnegative and diagonal.

- **EigvalueDecopoOnePass.m** - This is implementation of Algorithm 5.6 from paper [11].

  This function takes an Hermitian matrix $A$, a random test matrix $\Omega$, a sample matrix $Y = A\Omega$, and an orthonormal matrix $Q$ that can be generated by FixedRank.m, this computes an approximate eigenvalue decomposition $A \approx U\Lambda U^*$. This algorithm requires only one pass for the input matrix $A$; Comparing to previous algorithms it takes less storage.

- **BasicRandSVD.m** - This in implementation of RSVD from paper [13].

  This function takes an $m \times n$ matrix $A$, a target rank $k$, and an over-sampling parameter $p$ and computes matrices $U$, $D$, and $V$ in an approximate rank-$(k+p)$ SVD of $A$ (so that $U$ and $V$ are orthonormal, $D$ is diagonal, and $A \approx UDV^*$.

  To generate an input hermitian matrix, please visit and download the [random hermitian matrix generator function](#)[12].

- **AERandSVD.m** - This is implementation of ALGORITHM: ACCURACY ENHANCED RANDOMIZED SVD from paper [13].

  This function takes an $m \times n$ matrix $A$, a target rank $k$, an over-sampling parameter $p$, and a exponent $q$. It computes matrices $U$, $D$, and $V$ in an approximate rank-$(k+p)$ SVD of $A$ (so that $U$ and $V$ are orthonormal, $D$ is diagonal, and $A \approx UDV^*$. This algorithm is more accurate Compares to **BasicRandSVD.m**.

- **AEORandSVD.m** - This is implementation of ALGORITHM: ACCURACY ENHANCED RANDOMIZED SVD (WITH ORTHONORMALIZATION) from paper [13].

  This function takes an $m \times n$ matrix $A$, a target rank $k$, an over-sampling parameter $p$, and a exponent $q$. It computes matrices $U$, $D$, and $V$ in an approximate rank-$(k+p)$ SVD of $A$ (so that $U$ and $V$ are orthonormal, $D$ is diagonal, and $A \approx UDV^*$. This algorithm reduces the truncating error from the power iteration with large $q$ in **AERandSVD.m**.

- **SPRandEVDH.m** - This is implementation of ALGORITHM: SINGLE-PASS RANDOMIZED EVD FOR A HERMITIAN MATRIX from paper [13].

  This function takes an $n \times n$ Hermitian matrix $A$, a target rank $k$, and an over-sampling parameter $p$ and computes matrices $U$ and $D$ in an approximate rank-$k$ EVD of $A$ (so that $U$ is an orthonormal matrix, $D$ is a diagonal matrix, and $A \approx UDU^*$. This algorithm requires only one pass or the input matrix $A$; Comparing to previous algorithms it takes less storage.

  To generate an input hermitian matrix, please visit and download the [random hermitian matrix generator function](#)[12].

- **SPRandSVD.m** - This is implementation of ALGORITHM: SINGLE-PASS RANDOMIZED SVD FOR A GENERAL MATRIX from paper [13].

  This function takes an $m \times n$ matrix $A$, a target rank $k$, and an over-sampling parameter $p$ and computes matrices $U$, $D$, and $V$ in an approximate rank-$k$ SVD of $A$ (so that $U$ and $V$ are orthonormal, $D$ is diagonal, and $A \approx UDV^*$. This function extends **SPRandEVDH.m** to general input matrix.

- **randPowerMethod.m** - This is implementation of Algorithm 4 from paper [14].

  This function takes a Hermitian matrix $A$, a number $q$ for maximum number of iterations and a stopping tolerance $\epsilon$, it computes estimated $\xi$ for a maximum eigenvalue of $A$.

To generate an input hermitian matrix, please visit and download the random hermitian matrix generator function[12].

- **randomizedLanczos.m** - This is implementation of Algorithm 5 from paper [14].

    This function takes a Hermitian matrix $A$, a number $q$ for maximum number of iterations and computes estimated $(\xi; y)$ for a maximum eigenpair of $A$.

    To generate an input hermitian matrix, please visit and download the random hermitian matrix generator function[12].

| Function name | Reference | Algorithm # or name listed in the reference |
|---|---|---|
| EigvalueDecopoRow | [7] | 5.4 |
| EigvalueDecopoOnePass | [7] | 5.6 |
| EigvalueDecopoNystrom | [7] | 5.5 |
| DirectEigvalueDecopo | [7] | 5.3 |
| FastRandRF | [7] | 4.5 |
| Adaptive_randRF | [7] | 4.2 |
| AEORandSVD | [13] | ACCURACY ENHANCED RANDOMIZED SVD (WITH ORTHONORMALIZATION) |
| AERandSVD | [13] | ACCURACY ENHANCED RANDOMIZED SVD |
| BasicRandSVD | [13] | RSVD |
| FixedRank | [7] | Proto-Algorithm: Solving the Fixed -Rank Problem |
| randomizedLanczos | [14] | 5 |
| randPI | [7] | 4.3 |
| randPowerMethod | [14] | 4 |
| randRF | [7] | 4.1 |
| randSI | [7] | 4.4 |
| RandSVD | [7] | Prototype for Randomized SVD |
| SPRandEVDH | [13] | SINGLE-PASS RANDOMIZED EVD FOR A HERMITIAN MATRIX |
| SPRandSVD | [13] | SINGLE-PASS RANDOMIZED SVD FOR A GENERAL MATRIX |

**Sub Functions:**

- **SRFT.m** - This function proceed the subsampled random Fourier transform and produce an $n \times l$ matrix $\Omega$ in the form
$$\Omega = \sqrt{\frac{n}{l}} DFR,$$
  where
    * $D$ is an $n \times n$ diagonal matrix whose entries are independent random variables uniformly distributed on the complex unit circle.
    * $F$ is the $n \times n$ unitary discrete Fourier transform, whose entries take the values $f_{pq} = n^{-\frac{1}{2}} e^{-2\pi i (p-1)(q-1)/n}$ for $p, q = 1, 2, ..., n$.
    * $R$ is an $n \times l$ matrix that samples $l$ coordinates from $n$ uniformly at random.

- **SPD.m** - This function randomly generates a symmetric positive defined matrix $A$ with dimension $n \times n$.

- **maxeig.m** - This function returns the maximum absolute eigenvalue of input matrix $A$.

| Function name | For |
|---|---|
| SRFT | subsampled random Fourier transform |
| SPD | generates a random spd matrix |
| maxeig | returns the max eigenvalue |

**Drivers:**

- **imagedriver.m** - Image reconstruction using different methods of randomized SVD. This driver produce all experiment results for experiment 4.1 *Illustration of various Randomized SVD* listed in this report.

- **driver_bound.m** - Analysis on error bounds of different methods of the 1st stage of randomized SVD. This driver produce all experiment results for experiment 4.2 *Illustration of Error Bounds* listed in this report.

To compile and run this test driver, please download the following file from GitHub:

- controlledgap.m - This function takes input: $m, n$ as the size of desired matrix A, $r$ as the position of the gap and $gap$ for the size of the gap, and returns the testing matrix A that contains such a gap between its singular values at the defined position. Introduced in paper [15].

- **driver_sin.m** - Canonical angles for different methods of randomized SVD. This driver produce all experiment results for experiment 4.3 *Illustration of Canonical angles* listed in this report.

To compile and run this test driver, please download the following files from GitHub:

- controlledgap.m - This function takes input: $m, n$ as the size of desired matrix $A$, $r$ as the position of the gap and $gap$ for the size of the gap, and returns the testing matrix $A$ that contains such a gap between its singular values at the defined position. Introduced in paper [15].

- angle_bounds.m - This function takes the right singular vectors $V$, the starting guess $\Omega$, the singular values $s$, a target rank $k$ and the number of subspaces $q$, outputs both the bounds for $\sin(\theta(U_1, U_h))$ and $\sin(\theta(V_1, V_h))$. Introduced in paper [15].

- subspace_angles.m - This function computes the canonical angles between two subspaces $U$ and $U_h$. Introduced in paper[15].

| Driver file name | For |
|---|---|
| imagedriver | image experiment on chapter 3 |
| driver_bound | bound experiment on chapter 3 |
| driver_sin | canonical angles experiment on chapter 3 |

| Driver file name | For |
|---|---|
| imagedriver | image experiment on chapter 3 |
| driver_bound | bound experiment on chapter 3 |
| driver_sin | canonical angles experiment on chapter 3 |

**Data files:**

- **Sunflower.txt** - An 804 × 1092 matrix converted from a photo of Kansas Sunflowers.

| File name | For |
|---|---|
| Sunflower | matrix converted from a photo of Kansas Sunflowers |

# References

[1] Sanjoy Dasgupta and Anupam Gupta. An elementary proof of a theorem of Johnson and Lindenstrauss. *Random Structures & Algorithms*, 22(1):60–65, 2003.

[2] Petros Drineas, Ravi Kannan, and Michael W. Mahoney. Fast Monte Carlo algorithms for matrices II: Computing a low-rank approximation to a matrix. *SIAM Journal on computing*, 36(1):158–183, 2006.

[3] Petros Drineas and Michael W Mahoney. Lectures on randomized numerical linear algebra. *The Mathematics of Data*, 25:1, 2018.

[4] Petros Drineas, Michael W. Mahoney, and Shan Muthukrishnan. Sampling algorithms for $l_2$ regression and applications. In *Proceedings of the seventeenth annual ACM-SIAM symposium on Discrete algorithm*, pages 1127–1136, 2006.

[5] Petros Drineas, Michael W. Mahoney, Shan Muthukrishnan, and Tamás Sarlós. Faster least squares approximation. *Numerische Mathematik*, 117(2):219–249, 2011.

[6] Nathan Halko, Per-Gunnar Martinsson, and Joel A. Tropp. Finding structure with randomness: Stochastic algorithms for constructing approximate matrix decompositions. 2009.

[7] Nathan Halko, Per-Gunnar Martinsson, and Joel A. Tropp. Finding structure with randomness: Probabilistic algorithms for constructing approximate matrix decompositions. *SIAM Review*, 53(2):217–288, 2011.

[8] William B. Johnson and Joram Lindenstrauss. Extensions of Lipschitz mappings into a Hilbert space. *Contemporary Mathematics*, 26(1):189–206, 1984.

[9] Ravindran Kannan and Santosh Vempala. Randomized algorithms in numerical linear algebra. *Acta Numerica*, 26:95, 2017.

[10] Edo Liberty, Franco Woolfe, Per-Gunnar Martinsson, Vladimir Rokhlin, and Mark Tygert. Randomized algorithms for the low-rank approximation of matrices. *Proceedings of the National Academy of Sciences*, 104(51):20167–20172, 2007.

[11] Michael W. Mahoney. Randomized algorithms for matrices and data. *arXiv preprint arXiv:1104.5557*, 2011.

[12] Marcus. *Random Hermitian Matrix Generator*. MATLAB Central File Exchange, Retrieved April 28, 2021.

[13] Per-Gunnar Martinsson. Randomized methods for matrix computations. *The Mathematics of Data*, 25:187–231, 2019.

[14] Per-Gunnar Martinsson and Joel A. Tropp. Randomized numerical linear algebra: Foundations & algorithms. *arXiv preprint arXiv:2002.01387*, 2020.

[15] Arvind K. Saibaba. Randomized subspace iteration: Analysis of canonical angles and unitarily invariant norms. *SIAM Journal on Matrix Analysis and Applications*, 40(1):23–48, 2019.

[16] Lloyd N. Trefethen. The definition of numerical analysis. Technical report, Cornell University, 1992.

[17] John Von Neumann and Herman H. Goldstine. Numerical inverting of matrices of high order. *Bulletin of the American Mathematical Society*, 53(11):1021–1099, 1947.

[18] David P Woodruff. Sketching as a tool for numerical linear algebra. *arXiv preprint arXiv:1411.4357*, 2014.